\documentclass[12pt]{article}

\usepackage{graphicx} 
\usepackage{amsmath} 
\usepackage{amssymb,amsthm}
\usepackage{color}
\usepackage{tikz}
\usepackage{pgfplots}
\usetikzlibrary{math}

\newcommand{\R}{{\mathbb R}}

\newcommand{\C}{{\mathbb C}}

\newcommand{\hv}{\widehat{v}}
\newcommand{\hw}{\widehat{w}}
\newcommand{\hA}{\widehat{A}}
\newcommand{\hV}{\widehat{V}}
\newcommand{\hB}{\widehat{B}}
\newcommand{\hW}{\widehat{W}}
\newcommand{\hS}{\widehat{S}}
\newcommand{\hLambda}{\widehat{\Lambda}}

\let\epsilon\varepsilon
\let\theta\vartheta

\newtheorem{theorem}{Theorem}[section]\newtheorem{lemma}[theorem]{Lemma}

\newtheorem{corollary}[theorem]{Corollary}

\newtheorem{remark}[theorem]{Remark}

\title{On the Ginzburg-Landau approximation for quasilinear pattern forming reaction-diffusion-advection systems}
\author{Theo Belin, Guido Schneider \\
{\small
Institut f\"ur Analysis, Dynamik und Modellierung,} \\ {\small Universit\"at Stuttgart, Pfaffenwaldring 57, } \\ {\small 70569 Stuttgart, Germany}}

\date{January 22, 2026}

\begin{document}

\maketitle

\begin{abstract}
We prove that the Ginzburg-Landau equation correctly predicts the dynamics of quasilinear 
pattern-forming reaction-diffusion-advec\-tion systems, close to the 
first instability.
We present a simple theorem which is easily applicable for such systems and relies on key maximal regularity results.
The theorem is applied to the 
Gray-Scott-Klausmeier vegetation-water interaction model and its application to general reaction-diffusion-advection systems is discussed. 
\end{abstract}


\section{Introduction}


In the following we consider quasilinear 
pattern-forming reaction-diffusion-advection  systems
such as the
Gray-Scott-Klausmeier system which is given by
\begin{eqnarray}\label{gsk1}
\partial_t v & = & d \partial_x^2 v - bv + wv^2, \\
\partial_t w & = & \partial_x^2 (w^2) +c \partial_x w +a(1-w)-wv^2, 
\label{gsk2}
\end{eqnarray}
for $ v = v(x,t) \in \R $ describing the vegetation density, 
and  $ w = w(x,t) \in \R $ describing the water density,
for $ x \in \R $ and $ t \geq 0 $, 
with constants $ a , b \geq 0 $, $ c \in \mathbb{R} $, and $ d > 0 $.
The parameter  $ a $ roughly measures the amount of rainfall.
It is a model for the interaction between  vegetation and  water in semi-arid landscapes, describing, 
for instance, banded vegetation patterns on a sloped terrain.
The system  \eqref{gsk1}-\eqref{gsk2} shows quite rich dynamics and  recently attracted a lot of interest,
cf. \cite{rade1,rade2}.

In this paper, we are interested in the description of this system near the first instability in the parameter regime where a Turing instability occurs.
It is well known that, at least formally, the Ginzburg-Landau equation can be used for this purpose.
It appears as a universal amplitude equation that can be derived by a multiple scaling perturbation ansatz.
Error estimates which show that the Ginzburg-Landau 
equation makes correct predictions about the dynamics of the original 
pattern-forming system  
have been established  for semilinear systems in a number of papers,
starting with \cite{CE90,vH91,KSM92,Schn94a,Schn94b}.
Such estimates have been used, for instance, to prove the global existence 
of the solutions starting in a small neighborhood of
the weakly unstable origin, cf. \cite{Schn99JMPA}.
For an overview and an introduction to the theory see the textbook 
\cite{SU17book}. 

The only approximation result for quasilinear systems,
which is known to us,  can be found in \cite{Zi14PhD}. There, the Ginzburg-Landau approximation 
was justified for the Marangoni problem in Sobolev spaces $ H^n $.
The considered system \eqref{gsk1}-\eqref{gsk2} is also quasilinear due to the porous medium term $ \partial_x^2 (w^2) = 2w \partial_x^2 w + ( \partial_x w)^2 $ and is therefore parabolic only if $ w > 0 $.
 
This paper provides an easy-to-use approximation result for quasilinear pattern-forming reaction-diffusion-advection systems. We present our approach for the Gray-Scott-Klausmeier system and then formulate the general approximation result.

The plan of the paper is as follows. In Section \ref{sec2} we recall the 
linear stability analysis of the trivial solution 
for \eqref{gsk1}-\eqref{gsk2}
and recover Turing and Turing-Hopf 
instabilities. In Section \ref{sec3} we derive  
the Ginzburg-Landau equation associated to the Turing instability.
In Section \ref{secA} we introduce the functional analytic set-up which we use.
In Section \ref{sec4} we provide
estimates for the residual associated to the Ginzburg-Landau approximation, i.e.,
estimates for the terms which remain after inserting the approximation in the equations.
In Section \ref{sec5} we derive
the equations for the error made by this formal approximation. 
The final error estimates will follow in Section \ref{sec7}. 
The version for the general quasilinear pattern forming reaction-diffusion-advection systems can be found in Section \ref{secgeneral}.
The paper is closed with a discussion and outlook in Section \ref{sec8}.
In this last section we summarize our approach and outline the novelty of our method. 
\medskip

{\bf Acknowledgement.}
The work  is partially supported by the Deutsche Forschungsgemeinschaft
DFG through the cluster of excellence  'SimTech'  under EXC 2075-390740016.

\section{Linear stability analysis}
\label{sec2}

The time and space independent solutions $ w(x,t) = w $, $ v(x,t) = v \in \mathbb{R} $ 
of \eqref{gsk1}-\eqref{gsk2} satisfy 
$$ 
a(1-w)-wv^2 = 0 , \qquad - bv + wv^2 = 0 .
$$
We first find the trivial solution $ (v_0^\star,w_0^\star) = (0,1) $. Secondly, we find  $ vw = b $  and $ a(1-w) = b v $
which implies $ aw(1-w) = b^2 $ resp. $ a w^2 - a w + b^2 = 0 $. 
Therefore, we obtain the fixed points 
$$ 
w^\star_{\pm} = \frac{a \pm \sqrt{a^2 - 4 ab^2}}{2a} = \frac12 \pm \sqrt{\frac14 -  b^2/a} \qquad v_{\pm}^{\star} = b/w^\star_{\pm}
$$ 
for $ \frac14 -  b^2/a \geq 0 $ or equivalently $  a \geq 4 b^2 $, i.e., at $ a = 4 b^2 $ there is a saddle-node bifurcation.
The state $ (v_0^\star, w_0^\star) = (0,1) $ with zero vegetation represents the desert while the equilibria $ (w_+^\star,v_+^\star) $ and $ (w_-^\star,v_-^\star) $ represent co-existing homogeneously vegetated states.

The linearization around a fixed point $ (v,w) = (v^\star,w^\star) $ is given by 
\begin{eqnarray}\label{gsk1lin}
\partial_t v & = & d \partial_x^2 v - bv +  2 w^\star v^\star v + (v^\star)^2 w, \\
\partial_t w & = & 2 w^\star \partial_x^2 w +c \partial_x w - aw-  2 w^\star v^\star v - (v^\star)^2 w, 
\label{gsk2lin}
\end{eqnarray}
We immediately see that the desert state $ (v_0^\star,w_0^\star) = (0,1) $ is asymptotically stable in the PDE system for all $ a > 0 $ and 
$ b > 0 $. 
One of the bifurcating fixed points is stable and one is unstable in the ODE system.
To analyze the stability of the fixed points in the PDE system we consider 
the solutions 
$$ 
(v,w) = (\widehat{v},\widehat{w}) e^{ikx+\lambda t}
$$ 
of the linearized system \eqref{gsk1lin}-\eqref{gsk2lin}, where the eigenvalues $ \lambda = \lambda(k) $ are determined by the solutions of 
$$
 \det \left(\begin{array}{cc}
 -dk^2 -b+2w^\star v^\star -\lambda & (v^\star)^2  \\
  -2w^\star v^\star &  
-2w^\star k^2 +cik-a-(v^\star)^2 -\lambda  \\
\end{array}
\right) = 0 .
$$

Note that the $\lambda_{1, 2} : \R \mapsto \C$ satisfy $\lambda_i(-k) = \overline{\lambda_i(k)}$, $\Re(\lambda_1) \geq \Re(\lambda_2)$ and 
\begin{equation}
  \label{eq:asymptotics_lambda}
    \lambda_1(k) \sim  - dk^2,\qquad
    \lambda_2(k) \sim - 2 w^\star k^2,
    \end{equation}
for $ k \rightarrow \pm \infty $.

\begin{figure}[htbp] 
    \centering
    \tikzmath{
      \acrit = 0.2412; 
      \b = 0.2; 
      \c = 0;
      \d = 0.018;
      \ep = 0.02;
      function delta(\a) {
        return sqrt(1/4 - \b^2/\a);
      };
      function w(\a) {
        return 1/2 - delta(\a);
      };
      function v(\a){
        return \b/w(\a);
      };
      function det(\a){
        return 2*w(\a)*(v(\a))^3;
      };
      function pv(\k) {
        return \d * \k^2 - \b;
      };
      function pw(\a, \k) {
        return 2* w(\a) * \k^2 + \a + (v(\a))^2;
      };
    }
  
    \begin{tikzpicture}
  
      \begin{axis}[
          axis lines = middle,
          xlabel = $k$,
          ylabel = {$\lambda_1(k)$},
          legend pos = outer north east,
          legend style = {
            legend cell align = left,
          },
          domain=-4:4,
          samples=200,
      ]
  
  
      \addplot[red, thick] {-(pv(x) + pw(\acrit - \ep, x))/2 + sqrt((pv(x) - pw(\acrit - \ep, x))^2/4 - det(\acrit - \ep))};
      \addlegendentry{$a = a_{\text{crit}} - \epsilon^2$}
  
      \addplot[blue, very thick] {-(pv(x) + pw(\acrit, x))/2 + sqrt((pv(x) - pw(\acrit, x))^2/4 - det(\acrit))};
      \addlegendentry{$a = a_{\text{crit}}$}
  
      \addplot[green!30!black, thick] {-(pv(x) + pw(\acrit + \ep, x))/2 + sqrt((pv(x) - pw(\acrit + \ep, x))^2/4 - det(\acrit + \ep))};
      \addlegendentry{$a = a_{\text{crit}} + \epsilon^2$}
      
      \end{axis}
    \end{tikzpicture}

     \caption{For 
  $ a_{\text{crit}} = 0.2412 $, $ b = 0.2 $, $ c = 0 $, 
  $ d = 0.018 $ we have a spectral situation necessary for the derivation of a 
  Ginzburg-Landau equation. The figure shows the largest eigenvalue $ \lambda_1(k) \in \R$,
  as a function over the Fourier wave numbers.}
     \label{fig1}
  \end{figure}
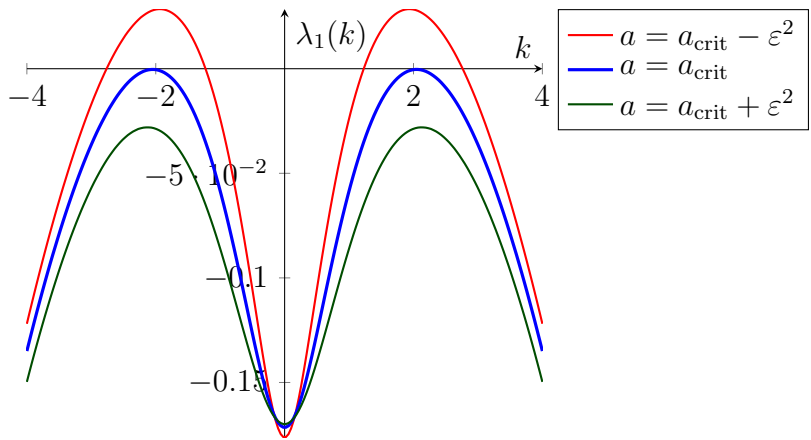

The fixed point $ (w_-^\star,v_-^\star) $ shows the kind of instability we are interested in.
In the following we
set $ (w^\star,v^\star) = (w_-^\star,v_-^\star) $, $\lambda_{1, 2} = \lambda^-_{1, 2}$ and 
 fix $ b $, $ c $, and $ d $, and take $ a $ as a bifurcation parameter.  
There is an open set of parameters $ b $, $ c $, and $ d $ where by decreasing $ a $ a Turing or a Turing-Hopf 
bifurcation occurs, cf. 
\cite[Figure 2.1]{rade2}. 
The wave number where this instability occurs is called $ k_{c} $ in the following.
Moreover, we set $ c = 0 $ and concentrate in the following analysis on the case of a pure subcritical Turing bifurcation, cf. Figure \ref{fig1}.  

\section{Derivation of the Ginzburg-Landau equation}
\label{sec3}

We consider the 
Gray-Scott-Klausmeier system \eqref{gsk1}-\eqref{gsk2} with $ (v,w) = (v^\star,w^\star) $ as new origin
and introduce the deviation 
$$ 
(\widetilde{v},\widetilde{w}) = (v,w) - (v^\star,w^\star) .
$$ 
The deviation satisfies 
\begin{equation} \label{OS}
\partial_t V = \Lambda V + B_2(V,V) + B_3(V,V,V),
\end{equation}
with $ V = (v,w)^T $,
$$ 
\Lambda V = \left(\begin{array}{c}
d \partial_x^2 v - bv +  2 w^\star v^\star v + (v^\star)^2 w \\ 
2 w^\star \partial_x^2 w +c \partial_x w - aw-  2 w^\star v^\star v - (v^\star)^2 w
\end{array} \right)
$$
and 
$$ 
B_2(V,V) = \left(\begin{array}{c}
   w^\star v^2 + 2 v^\star v w    \\ 
  \partial_x^2 (w^2) - w^\star v^2 - 2 v^\star v w  
\end{array} \right),
\qquad 
B_3(V,V,V) = \left(\begin{array}{c}
   w v^2   \\ 
 - w v^2 
\end{array} \right),
$$
where we skipped the \ $ \widetilde{}$ s for notational simplicity and where 
w.l.o.g. we assume the symmetry of $ B_2 $ and $ B_3 $ in their arguments.
For the derivation of the Ginzburg-Landau equation we consider 
the Fourier transformed system 
\begin{equation} \label{sysFT}
\partial_t \hV = \hLambda \hV + \hB_2(\hV,\hV) + \hB_3(\hV,\hV,\hV),
\end{equation}
with $ \hV = (\hv,\hw)^T $,
$$ 
\hLambda(k) \hV = \left(\begin{array}{c}
- d k^2 \hv - b\hv +  2 w^\star v^\star \hv + (v^\star)^2 \hw \\ 
- 2 w^\star k^2 \hw +c ik  \hw - a\hw-  2 w^\star v^\star \hv - (v^\star)^2 \hw
\end{array} \right)
$$
and
\begin{eqnarray*}
\hB_2(k)(\hV,\hV) & = & \left(\begin{array}{c}
   w^\star\hv^{*2} + 2 v^\star*\hv *\hw    \\ 
- k^2 (\hw^{*2}) - w^\star\hv^{*2} - 2 v^\star \hv *\hw  
\end{array} \right),
\\
\hB_3(\hV,\hV,\hV) & = & \left(\begin{array}{c}
   \hw *\hv^{*2}   \\ 
 - \hw *\hv^{*2} 
\end{array} \right).
\end{eqnarray*}
Next we 
introduce the small bifurcation parameter $ 0 < \varepsilon \ll 1 $
by
$$ 
a   = a_{crit} - \varepsilon^2.
$$ 
We can diagonalize 
$ \hLambda(k) $ near $ k_c $ by $ \hV(k) = \hS(k) \hW(k) $
such that  
$$
\partial_t \hW = \widehat{D}(k) \hW + \hS(k)^{-1} \hB_2(\hV,\hV) + \hS(k)^{-1} \hB_3(\hV,\hV,\hV),
$$
with 
$$  
\widehat{D}(k) = \left(\begin{array}{cc}
\lambda_1(k) & 0  \\  0 & \lambda_2(k)
\end{array} \right),
$$ 
where we recall that $ \lambda_1(k) \geq  \lambda_2(k) $ are the real eigenvalues  
of  $ \hLambda $. The associated eigenvectors are called $ f_1(k) $ and $ f_2(k) $.
Then we make the ansatz 
\begin{eqnarray*}
\varepsilon \widehat{\psi}_{GL}(k,t) & = & \varepsilon \varepsilon^{-1}  \hA_1\left(\frac{k-k_c}{\varepsilon},\varepsilon^2 t\right) f_1(k) + \varepsilon \varepsilon^{-1}  \hA_{-1}\left(\frac{k+k_c}{\varepsilon},\varepsilon^2 t\right) f_1(k)  \\ 
&&+  \varepsilon^2 \varepsilon^{-1}  \hV_2\left(\frac{k-2k_c}{\varepsilon},\varepsilon^2 t\right) +
\varepsilon^2 \varepsilon^{-1}  \hV_{-2}\left(\frac{k+2k_c}{\varepsilon},\varepsilon^2 t\right) \\
&& + \varepsilon^2 \varepsilon^{-1}  \hV_0\left(\frac{k}{\varepsilon},\varepsilon^2 t\right),
\end{eqnarray*}
where $A_{-1} = \overline{A_1}$, $V_{-2} = \overline{V_2}$, and $ V_0 = \overline{V_0} $. %

Plugging in this ansatz in \eqref{sysFT}, 
equating the coefficients at different powers of $ \varepsilon $ at different  
wave numbers to zero, and using the abbreviations 
$ K_j = \frac{k - j k_c}{\varepsilon} $ and $ T = \varepsilon^2 t $,
gives the approximation equations.
At $ k = k_c $ and $ \varepsilon^0 $ and $ \varepsilon^1 $ all coefficients 
vanish. At $ \varepsilon^2 $ and $ f_1 $ we find
for $ \hA_j = \hA(K_j,T) $, $\hV_j = \hV_j(K_j, T)$, that 
\begin{eqnarray*}
\partial_T \hA_1 & = &\alpha_0  \hA_1 -  \alpha_2 K_1^2  \hA_1  + 2 f_1^*(k_c) \cdot \hB_2(k_c)(\hA_1 f_1(k_c),\hV_0) 
\\ && + 2 f_1^*(k_c) \cdot \hB_2(k_c)(\hA_{-1} f_1(-k_{c}),\hV_2) \\ && + 3 f_1^*(k_c) \cdot \hB_2(\hA_{-1} f_1(-k_{c}), \hA_1 f_1(k_c),\hA_1 f_1(k_c)),
\end{eqnarray*}
with 
$$
\alpha_0 = \lim_{\varepsilon \to 0} \varepsilon^{-2} \partial_a \lambda_1|_{k = k_c, a= a_{crit} + \varepsilon^{2}}, \qquad 
\alpha_2 = - \partial_k^2 \lambda_1|_{k = k_c, a= a_{crit}}/2,
$$
where $ f_j^* $ is the adjoint eigenvector to $ f_j $. 
At $ k = 0 $ and $ \varepsilon^1 $ we find 
$$ 
0 =  \hLambda(0) \hV_0 + 2 \hB_2(0)(\hA_1 f_1(k_c), \hA_{-1} f_1(-k_{c})).
$$ 
At $ k = 2 k_c $ and $ \varepsilon^1 $ we find 
$$ 
0 =  \hLambda(2 k_c) \hV_2 +  \hB_2(2 k_c)(\hA_1 f_1(k_c), \hA_1 f_1(k_c)).
$$ 
The last two equations can be solved w.r.t. $ \hV_0 $ and $ \hV_2$. Inserting this in
the equation for $ \hA_1 $ yields the Ginzburg-Landau equation
\begin{equation}\label{GLF}
\partial_T \hA_1  = \alpha_0  \hA_1 -  \alpha_2 K_1^2  \hA_1 + \alpha_3 \hA_1 * \hA_1 * \hA_{-1},
\end{equation} 
with explicitly computable coefficient $ \alpha_3 $.
In physical space the Ginzburg-Landau equation is given by 
\begin{equation}\label{GL}
\partial_T A_1  = \alpha_0  A_1 +  \alpha_2 \partial_X^2  A_1 + \alpha_3 |A_1|^2 A_1 .
\end{equation} 

\section{The functional analytic set-up}

\label{secA}

We 
follow \cite{BK92,BK94}.
For each $ r \geq 0 $ we introduce the  spaces 
$$ 
L^{\infty}_r = \{ \widehat{u} \in L^{\infty}(\R,\C) :  \|\widehat{u} \|_{L^{\infty}_r} < \infty \},
$$ 
equipped  with the norm
$$
\|\widehat{u}\|_{L^{\infty}_r} = \sup_{k \in \R} |\widehat{u}(k) |(1+|k|^2)^{r/2}.
$$ 
\begin{lemma} \label{lemconv}
The spaces $ L^{\infty}_r(\R;\C) $ are closed under convolution if $ r > 1 $.
In detail, for all $ r > 1 $ there exists a $ C > 0 $ such that for all 
$ \widehat{u} , \widehat{v} \in L^{\infty}_{r} $ we have
$$ 
\|\widehat{u}*\widehat{v} \|_{L^{\infty}_{r}} \leq 
C \|\widehat{u} \|_{L^{\infty}_{r}}
\|\widehat{v} \|_{L^{\infty}_{r}}.
$$
\end{lemma}
\noindent
{\bf Proof.}
The statement follows from Young's inequality for convolutions 
$$ 
\|\widehat{u}*\widehat{v} \|_{L^{\infty}_{r}} \leq C 
(\|\widehat{u} \|_{L^1} \|\widehat{v} \|_{L^{\infty}_{r}} + 
\|\widehat{u}\|_{L^{\infty}_{r}} \|\widehat{v} \|_{L^1} )
$$ 
and 
Sobolev's embedding in Fourier space
$$ 
\|\widehat{u}\|_{L^1} \leq 
\| (1+|\cdot |^2)^{-r/2}\|_{L^1}
 \|\widehat{u} \|_{L^{\infty}_{r}}
\leq 
C \|\widehat{u}\|_{L^{\infty}_{r}}
$$ 
which holds since $ \| (1+|\cdot |^2)^{-r/2}\|_{L^1}  = C < \infty $ 
if $ r > 1 $. 
\qed
\medskip

For notational simplicity we prefer to work in physical space
 and so we define 
 $$ 
 \mathcal{X}^r = \{ u: \R \to \C : \widehat{u} \in L^{\infty}_r \}
 $$ 
 equipped with the norm 
 $$ 
 \| u \|_{\mathcal{X}^r} = \|\widehat{u} \|_{L^{\infty}_r}.
 $$ 
 A reformulation of Lemma \ref{lemconv} in terms of $ \mathcal{X}^r $ yields
 \begin{corollary} \label{collconv}
The spaces $ \mathcal{X}^r $ are closed under pointwise multiplication if $ r > 1 $.
In detail, for all $ r > 1 $ there exists a $ C_r > 0 $ such that for all $u, v \in \mathcal{X}^r$ we have
$$ 
\|u \cdot v \|_{\mathcal{X}^r} \leq 
C_r \|u \|_{\mathcal{X}^r}
\|v\|_{\mathcal{X}^r}.
$$
\end{corollary}

There is local existence and uniqueness for the Ginzburg-Landau equation 
\eqref{GL}
in these spaces.
\begin{lemma} \label{lemexist}
For all $ r > 1 $, $C_{GL} > 0$ there exists a $ T_0 = T_0(C_{GL}) > 0 $ such that the following holds.
Let $ A_0 \in \mathcal{X}^r $ with $ \| A_0 \|_{\mathcal{X}^r} \leq C_{GL} $.
Then there exists a unique solution $ A \in C([0,T_0],\mathcal{X}^r) $ 
of the Ginzburg-Landau equation 
\eqref{GL} with $ A(\cdot,0) = A_0(\cdot) $.
\end{lemma}
\noindent
{\bf Proof.} 
The proof follows the usual approach for semilinear parabolic equations, cf. 
\cite{He81}.
We apply the variation of constant formula to 
the Ginzburg-Landau equation. Then the right hand side is a contraction
in a ball of $ \mathcal{X}^r $
with radius $ 2 C_{GL} $ for $ T_0 > 0 $ sufficiently small.
\qed
\medskip

To see that the Ginzburg-Landau equation makes correct predictions
about the Gray-Scott-Klausmeier system we have to prove 
error estimates for all $ t \in [0,T_0/\varepsilon^2] $.
A simple application of  Gronwall's inequality  is not sufficient 
to obtain such estimates  due to the quadratic terms present 
in  the Gray-Scott-Klausmeier system.
We follow the existing literature 
and introduce so called  mode-filters to separate the critical 
modes from the exponentially damped modes.
Here, they can be   
defined as indicator functions in Fourier space by
\[
 \widehat{E}_j \widehat{u}= \left\{ \begin{array}{cl} 1 ,&  |k - j  k_c|\leq k_c/10,
 \\
 0 ,&  |k -j  k_c | >  k_c/10,
 \end{array}\right.
 \]
 by $  \widehat{E}_c =  \widehat{E}_1 +  \widehat{E}_{-1} $, 
and by
 $  \widehat{E}_s=1- \widehat{E}_{1}- \widehat{E}_{-1}$.
 We obviously have 
 \begin{lemma}
 The mode-filters  are smooth mappings in each $ L^{\infty}_{r} $ for $ r \geq 0 $. In detail, 
 for all $ r \geq 0 $ there exists a $ C  > 0 $ such that 
 $$ 
 \| \widehat{E}_{\pm 1}  \widehat{u} \|_{L^{\infty}_{r}} +  \| \widehat{E}_{c}  \widehat{u} \|_{L^{\infty}_{r}} 
 +  \| \widehat{E}_{s}  \widehat{u} \|_{L^{\infty}_{r}} 
 \leq C \|  \widehat{u}  \|_{L^{\infty}_{r}} .
 $$
 \end{lemma}
 A reformulation of Lemma \ref{lemconv} in  $ \mathcal{X}^r $-spaces yields
\begin{corollary} 
 The mode-filters are smooth mappings in each $ \mathcal{X}^r $ for $ r \geq 0 $. In detail, 
 for all $ r \geq 0 $ there exists a $ C  > 0 $ such that 
 $$ 
 \| E_{\pm 1}  u \|_{\mathcal{X}^r} +  \| E_{c}  u \|_{\mathcal{X}^r} 
 +  \|E_{s}  u \|_{\mathcal{X}^r} 
 \leq C \|  u  \|_{\mathcal{X}^r} ,
 $$
where $ E_j = \mathcal{F}^{-1} \widehat{E}_{j} \mathcal{F} $ for $ j = \pm 1,c,s $.
\end{corollary}

\section{Estimates for the residual}
\label{sec4}

The residual 
\begin{equation} \label{res}
\textrm{Res}(V) = \left(\begin{array}{c}
  \textrm{Res}_v(V)  \\ 
\textrm{Res}_w(V)
\end{array} \right) = - \partial_t V + \Lambda V + B_2(V,V) + B_3(V,V,V)
\end{equation}
contains all terms which do not cancel after inserting the 
approximation into  \eqref{OS}. It counts how much an approximation fails to satisfy the equations \eqref{OS}.
If $ \textrm{Res}(V) = 0 $ then $ V $ is a solution of \eqref{OS}. 

At this point our choice of spaces leads to some additional work
due to the scaling properties of the Fourier transform.
By the construction of the Ginzburg-Landau approximation in Section 
\ref{sec3} all terms of
order $ \mathcal{O}(\varepsilon^2) $ have canceled around $ \pm k_c $ 
and all terms of
order $ \mathcal{O}(\varepsilon) $ have canceled at all other places.
Due to our choice of norm this is also true for the $ E_c $-part and the 
$ E_s $-part, respectively.
However, we also need to cancel the terms of order $ \mathcal{O}(\varepsilon^3) $ and $ \mathcal{O}(\varepsilon^2) $, respectively. This can be achieved, as usual, by adding higher order terms 
to the approximation, cf. \cite[Section 10.5]{SU17book}.
We choose the improved approximation
\begin{eqnarray*}
\varepsilon \widehat{\psi}(k,t) & = & 
\varepsilon \widehat{\psi}_{GL}(k,t) \\ && + 
\varepsilon^2 \varepsilon^{-1}  \hA_{1,1}\left(\frac{k-k_c}{\varepsilon},\varepsilon^2 t\right) f_1(k) + \varepsilon^2 \varepsilon^{-1}  \hA_{-1,1}\left(\frac{k+k_c}{\varepsilon},\varepsilon^2 t\right) f_1(k)  \\ &&
+ 
\varepsilon^3 \varepsilon^{-1}  \hA_{1, 2}\left(\frac{k-k_c}{\varepsilon},\varepsilon^2 t\right) f_2(k) + \varepsilon^3 \varepsilon^{-1}  \hA_{-1, 2}\left(\frac{k+k_c}{\varepsilon},\varepsilon^2 t\right) f_2(k)  \\
&&+  \varepsilon^3 \varepsilon^{-1}  \hV_{2,1}\left(\frac{k-2k_c}{\varepsilon},\varepsilon^2 t\right) +
\varepsilon^3 \varepsilon^{-1}  \hV_{-2,1}\left(\frac{k+2k_c}{\varepsilon},\varepsilon^2 t\right) \\
&& + \varepsilon^3 \varepsilon^{-1}  \hV_{0,1}\left(\frac{k}{\varepsilon},\varepsilon^2 t\right) 
\\ &&+  \varepsilon^3 \varepsilon^{-1}  \hV_{3}\left(\frac{k-3k_c}{\varepsilon},\varepsilon^2 t\right) +
\varepsilon^3 \varepsilon^{-1}  \hV_{-3}\left(\frac{k+3k_c}{\varepsilon},\varepsilon^2 t\right),
\end{eqnarray*}
where the $ A_{\pm 1,1} $ satisfy linearized Ginzburg-Landau equations and 
where the other functions are determined by linear inhomogeneous equations, cf. \cite[Section 10.5]{SU17book}.
We have the following estimates 
for the improved Ginzburg-Landau  approximation.
\begin{lemma} \label{lem51}
For all $ C_{GL} > 0 $, $ T_0 > 0 $, $ r > 1 $ there exist $ \varepsilon_0 > 0 $ 
and $ C_2 > 0 $ such that the following holds.
Let $ A \in C([0,T_0],\mathcal{X}^{r+6}) $ be a solution of the 
Ginzburg-Landau equation \eqref{GL} satisfying 
\begin{equation} \label{eq19}
\sup_{T \in [0,T_0/\varepsilon^2]} \| A(\cdot,T) \|_{\mathcal{X}^{r+6}}
\leq C_{GL}.
\end{equation}
Then  for all $ \varepsilon \in (0,\varepsilon_0)
$ we have 
$$
\sup_{t \in [0,T_0/\varepsilon^2]} \| E_c \textrm{Res}(\varepsilon \psi) \|_{\mathcal{X}^r} \leq C_2 \varepsilon^4 
$$ 
and 
$$
\sup_{t \in [0,T_0/\varepsilon^2]} \| E_s \textrm{Res}(\varepsilon \psi) \|_{\mathcal{X}^r} \leq C_2 \varepsilon^3 .
$$ 
\end{lemma}
\noindent
{\bf Proof.}
The proof follows almost line for line the one given in \cite[Section 10.5]{SU17book} for reaction-diffusion systems.
\qed 

\section{The equations for the error}
\label{sec5} \label{sec6}

We start with \eqref{OS}. We write the solution $ V $ as a sum of the improved approximation $ \varepsilon \psi $ 
and an error $  \varepsilon^2 R $. We split the approximation as
$$ 
\varepsilon \psi = \varepsilon \psi_c + \varepsilon^2 \psi_s,
$$ 
where $\psi_c = E_c \psi$ and $\varepsilon \psi_s = E_s \psi$.  Thus, the support of $ \psi_c $ in Fourier space is contained in the support of $ \widehat{E}_c $ 
and that of $ \psi_s $ is contained in the support of $ \widehat{E}_s $.
Moreover, we split the error $\varepsilon^2 R = V - \varepsilon \psi$, as
$$ 
\varepsilon^2 R = \varepsilon^2 R_c + \varepsilon^3 R_s,
$$ 
where $ R_c = E_c R $ and $ \varepsilon R_s = E_s R $.  Thus, the support of $R_c$ in Fourier space is contained in the support of $ \widehat{E}_c $ 
and that of $ R_s $ is contained in the support of $ \widehat{E}_s $.
We note that $ R_c $ and $ R_s $ are solutions of 
\begin{eqnarray} \label{erreq1}
\partial_t R_c & = & \Lambda R_c + L_{2,c}(\varepsilon \psi,\varepsilon^2 R)
+ N_{2,c}(\varepsilon \psi,\varepsilon^2 R)  
 \\ && \qquad + \varepsilon^{-2} E_c N_3(\varepsilon \psi,\varepsilon^2 R) + \varepsilon^{-2} E_c \textrm{Res}(\varepsilon \psi), 
 \nonumber
 \\ \label{erreq2}
\partial_t R_s & = & \Lambda R_s + L_{2,s}(\varepsilon \psi,\varepsilon^2 R)
+ N_{2,s}(\varepsilon \psi,\varepsilon^2 R)   \\ && \qquad  + \varepsilon^{-3} E_s N_3(\varepsilon \psi,\varepsilon^2 R)
+ \varepsilon^{-3} E_s \textrm{Res}(\varepsilon \psi), \nonumber
\end{eqnarray}
where 
\begin{eqnarray*}
L_{2,c}(\varepsilon \psi,\varepsilon^2 R) & = & 2 \varepsilon^2  E_c B_2(\psi_c, R_s) 
+ 2 \varepsilon^2  E_c B_2(\psi_s, R_c)  + 2 \varepsilon^3  E_c B_2(\psi_s, R_s) , \\ 
N_{2,c}(\varepsilon \psi,\varepsilon^2 R) & = & 2 \varepsilon^3  E_c B_2(R_c, R_s) + \varepsilon^4  E_c B_2(R_s, R_s), \\
L_{2,s}(\varepsilon \psi,\varepsilon^2 R) & = & 2   E_s B_2(\psi_c, R_c) + 2 \varepsilon  E_s B_2(\psi_c, R_s)   \\ && \qquad  + 2 \varepsilon  E_s B_2(\psi_s, R_c)  + 2 \varepsilon^2  E_s B_2(\psi_s, R_s) ,\\
N_{2,s}(\varepsilon \psi,\varepsilon^2 R) & = &  \varepsilon  E_s B_2(R_c, R_c) + 2 \varepsilon^2  E_s B_2(R_c, R_s) + \varepsilon^3  E_s B_2(R_s, R_s),\\
N_3(\varepsilon \psi,\varepsilon^2 R) & = & 3 B_3(\varepsilon \psi ,\varepsilon \psi ,\varepsilon^2 R) 
+ 3 B_3(\varepsilon \psi ,\varepsilon^2 R ,\varepsilon^2 R) + B_3(\varepsilon^2 R ,\varepsilon^2 R ,\varepsilon^2 R) .
\end{eqnarray*}
Note that $ E_c B_{2}(\psi_c,R_c) =  E_c B_{2}(R_c,R_c) = 0 $ due to disjoint supports 
of $ E_c $ on the one hand and $ B_{2}(\psi_c,R_c)  $ and $ B_{2}(R_c,R_c)  $ on the other hand 
in Fourier space.

We have the estimates for the cubic terms 
\begin{eqnarray*}
\| \varepsilon^{-2} E_c N_3(\varepsilon \psi,\varepsilon^2 R)  \|_{\mathcal{X}^r} & \leq & C \varepsilon^2 \| R \|_{\mathcal{X}^r}  + C  \varepsilon^3 \| R \|_{\mathcal{X}^r}^2 + C  \varepsilon^4 \| R \|_{\mathcal{X}^r}^3, \\  
\| \varepsilon^{-3} E_s N_3(\varepsilon \psi,\varepsilon^2 R)  \|_{\mathcal{X}^r} & \leq & C \varepsilon \| R \|_{\mathcal{X}^r}  + C  \varepsilon^2 \| R \|_{\mathcal{X}^r}^2 + C  \varepsilon^3 \| R \|_{\mathcal{X}^r}^3,
\end{eqnarray*}
for all $ r > 1 $. 
Since $ E_c $ has compact support in Fourier space it is smooth from $ \mathcal{X}^0 $ to every $ \mathcal{X}^r $, $r \geq 0$.
Therefore, we easily find 
\begin{eqnarray*}
 \| L_{2,c}(\varepsilon \psi,\varepsilon^2 R)
+ N_{2,c}(\varepsilon \psi,\varepsilon^2 R)   \|_{\mathcal{X}^r} 
 \leq  C \varepsilon^2 \| R \|_{\mathcal{X}^r}  + C  \varepsilon^3 \| R \|_{\mathcal{X}^r}^2 ,
\end{eqnarray*}
for all $ r > 1 $.
The semilinear part in $ L_{2,s} $ and $ N_{2,s} $ can be handled similarly. Therefore, we split
$ B_2 $ in a semilinear part $ B_{2,s} $ and in a quasilinear part $ B_{2,q} $, with  
$$ 
B_{2,s}(V,V) = \left(\begin{array}{c}
   w^\star v^2 + 2 v^\star v w    \\ 
 - w^\star v^2 - 2 v^\star v w  
\end{array} \right), \qquad 
B_{2,q}(V,V) = \left(\begin{array}{c}
  0   \\ 
  \partial_x^2 (w^2)  
\end{array} \right).
$$
The associated semilinear parts $ L_{2,s,s} $ and $ N_{2,s,s} $ can be estimated easily as
\begin{eqnarray*}
 \| L_{2,s,s}(\varepsilon \psi,\varepsilon^2 R)
+ N_{2,s,s}(\varepsilon \psi,\varepsilon^2 R)   \|_{\mathcal{X}^r} 
 \leq  C  \| R_c \|_{\mathcal{X}^r}  + C \varepsilon \| R_s \|_{\mathcal{X}^r}  + C  \varepsilon \| R \|_{\mathcal{X}^r}^2 
\end{eqnarray*}
for all $ r > 1 $.
The associated quasiilinear parts $ L_{2,s,q} $ and $ N_{2,s,q} $ can be estimated easily as
\begin{eqnarray*}
 \| L_{2,s,q}(\varepsilon \psi,\varepsilon^2 R)
+ N_{2,s,q}(\varepsilon \psi,\varepsilon^2 R)   \|_{\mathcal{X}^{r-2}} 
 \leq  C  \| R_c \|_{\mathcal{X}^r}  + C \varepsilon \| R_s \|_{\mathcal{X}^r}  + C  \varepsilon \| R \|_{\mathcal{X}^r}^2 
\end{eqnarray*}
for all $ r \geq 2 $.
\begin{remark}{\rm
We used 
$$ 
\|R_{j_1}(\cdot,t) R_{j_2}(\cdot,t) \|_{\mathcal{X}^r} \leq 
C \|R_{j_1}(\cdot,t)  \|_{\mathcal{X}^r}
\|R_{j_2}(\cdot,t)  \|_{\mathcal{X}^r},
$$
from Corollary \ref{collconv} but 
$$ 
\|\psi_{j_1}(\cdot,t) R_{j_2}(\cdot,t) \|_{\mathcal{X}^r} \leq 
C \|\widehat{\psi}_{j_1}(\cdot,t)  \|_{L^1_r}
\|R_{j_2}(\cdot,t)  \|_{\mathcal{X}^r},
$$
since 
\begin{equation}\label{embedding}
\|\widehat{\psi}_{j_1}(\cdot,t)  \|_{L^1_r} 
=\int |\widehat{\psi}_{j_1}(k,t) | (1+k^2)^{r/2}  dk = \mathcal{O}(1)
\end{equation}
in contrast to  $ \|\psi_{j_1}(\cdot,t)  \|_{\mathcal{X}^r} = \mathcal{O}(1/\varepsilon) $.
Estimate \eqref{embedding} is a direct consequence 
of $ \| \widehat{A} \|_{L^1_r} \leq C  \| \widehat{A} \|_{L^{\infty}_{r_A}} $
if $ r_A- r > 1 $.
}\end{remark}
From Lemma \ref{lem51} we additionally have 
$$
\sup_{t \in [0,T_0/\varepsilon^2]} \| \varepsilon^{-2} E_c \textrm{Res}(\varepsilon \psi) \|_{\mathcal{X}^r} \leq C \varepsilon^2 
\quad
\textrm{and} 
\quad
\sup_{t \in [0,T_0/\varepsilon^2]} \| \varepsilon^{-3} E_s \textrm{Res}(\varepsilon \psi) \|_{\mathcal{X}^r} \leq C  .
$$ 

\section{The error estimates}
\label{sec7}

For solving the equations of the error we define the spaces 
$$ 
 \mathcal{Y}^{r,\eta} = \{ u \in C([0,T_0/\varepsilon^2],\mathcal{X}^r)
:
 \| {u} \|_{\mathcal{Y}^{r,\eta}} < \infty \}
 $$ 
 equipped with the norm 
 $$ 
 \| u \|_{ \mathcal{Y}^{r,\eta}} = \sup_{t \in [0,T_0/\varepsilon^2]} e^{-\eta t}\|u(\cdot,t) \|_{\mathcal{X}^r} = \sup_{t \in [0,T_0/\varepsilon^2]} e^{-\eta t}\|\widehat{u}(\cdot,t) \|_{L^{\infty}_r}.
 $$ 
 In the following we choose $ \eta = \widetilde{\eta} \varepsilon^2 $
 with $ \widetilde{\eta} = \mathcal{O}(1) > 1 $ fixed, but sufficiently large. 
 Since $ e^{\widetilde{\eta} \varepsilon^2 t} = \mathcal{O}(1) $ for 
 $ t \in [0,T_0/\varepsilon^2]$ all estimates from the last section transfer 
 one-to-one from $ \mathcal{X}^r $-spaces to $ \mathcal{Y}^{r,\widetilde{\eta} \varepsilon^2} $-spaces. For example we have 
  \begin{corollary} \label{collconvy}
The spaces $ \mathcal{Y}^{r,\widetilde{\eta} \varepsilon^2}  $ are closed under pointwise multiplication if $ r > 1 $.
In detail, for all  $ r > 1 $ there exist $ \varepsilon_0 > 0 $ and  $ C > 0 $ such that for all $ \varepsilon \in (0,\varepsilon_0) $ and $  \widetilde{\eta} > 0 $, we have  
$$ 
\|u \cdot v \|_{\mathcal{Y}^{r,\widetilde{\eta} \varepsilon^2} } \leq 
C e^{\widetilde{\eta}}\|u \|_{\mathcal{Y}^{r,\widetilde{\eta} \varepsilon^2} }
\|v \|_{\mathcal{Y}^{r,\widetilde{\eta} \varepsilon^2} }.
$$
\end{corollary}
Next we estimate the linear operator 
$$
\mathcal{K}_c :\left\{ \begin{array}{c} f \mapsto R_c, \\
\mathcal{Y}^{r,\widetilde{\eta} \varepsilon^2} \to \mathcal{Y}^{r,\widetilde{\eta} \varepsilon^2},
\end{array} \right.
$$  
where $ R_c $  is defined through 
\begin{eqnarray*}
\partial_t R_c & = & \Lambda R_c 
+  \varepsilon^2 E_c f  
\end{eqnarray*}
and the linear operator 
$$
\mathcal{K}_s :\left\{ \begin{array}{c} f \mapsto R_s, \\
\mathcal{Y}^{r-2,\widetilde{\eta} \varepsilon^2} \to \mathcal{Y}^{r,\widetilde{\eta} \varepsilon^2},
\end{array} \right.
$$
where $ R_s $  is defined through 
\begin{eqnarray*}
\partial_t R_s & = & \Lambda R_s 
+ E_s f.
\end{eqnarray*}

Note it is straightforward to check that, for $r > 0$, the linear operator $\Lambda$ is closed, has dense domain in $\mathcal{X}^r$ and generates a strongly continuous semi-group $(e^{\Lambda t})_{t\geq 0}$ on $\mathcal{X}^r$. As a result the linear operators $\mathcal{K}_c$ and $\mathcal{K}_s$ are well-defined. For $ \mathcal{K}_c $ we prove 
\begin{lemma}
  \label{lem72}
For all  $ r \geq 0  $ there exist $ \varepsilon_0 > 0 $ and  $ C > 0 $ such that for all $ \varepsilon \in (0,\varepsilon_0) $ and all $ \widetilde{\eta} \geq 2 $ we have  
$$
\| \mathcal{K}_c f \|_{\mathcal{Y}^{r,\widetilde{\eta} \varepsilon^2} } \leq C (\widetilde{\eta}-1)^{-1}  \|  f \|_{
\mathcal{Y}^{0,\widetilde{\eta} \varepsilon^2}}.
$$
\end{lemma}
\noindent
{\bf Proof.} 
We set 
$$ 
W_c(\cdot,t) = e^{-\widetilde{\eta} \varepsilon^2 t} R_c(\cdot,t) 
$$
and find 
$$ 
\partial_t W_c(\cdot,t)  =  (\Lambda- \widetilde{\eta} \varepsilon^2) 
W_c(\cdot,t) 
+  \varepsilon^2  e^{-\widetilde{\eta} \varepsilon^2 t}  (E_c f) (\cdot,t). 
$$ 
Applying the variation of constant formula yields
\begin{eqnarray*}
W_c(t) & = & \varepsilon^2 \int_0^t e^{ (\Lambda- \widetilde{\eta} \varepsilon^2)  (t-\tau)} e^{-\widetilde{\eta} \varepsilon^2 \tau} E_c f(\tau) d\tau, 
\end{eqnarray*}
where the integral above is well-defined by strong continuity of the semi-group generated by $\Lambda$ and continuity of $f: [0, T_0/\varepsilon^2] \rightarrow \mathcal{X}^r$. \\

For the semigroup $ (e^{\Lambda t})_{t \geq 0} $ applied 
to the $ E_c $-part we find 
$$
\sup_{|k \pm k_c| \leq k_c/10} |  e^{(\Lambda- \widetilde{\eta} \varepsilon^2)t} | \leq C_{\Lambda} e^{(1-\widetilde{\eta} )\varepsilon^2 t} 
$$
%
with a constant $ C_{\Lambda} = \mathcal{O}(1) $.
Then for  each $ r \geq 0 $ 
we estimate   with $ f_c = E_c f $ that 
\begin{eqnarray*}
\lefteqn{\sup_{t \in [0,T_0/\varepsilon^2]}\|\widehat{W}_c(\cdot,t)  \|_{L^{\infty}_r} }\\ & = & 
\varepsilon^2 \sup_{t \in [0,T_0/\varepsilon^2]}\sup_{|k \pm k_c| \leq k_c/10}\left( \left| \int_0^t e^{(\Lambda- \widetilde{\eta} \varepsilon^2)(t-\tau)} e^{-\widetilde{\eta} \varepsilon^2 \tau}\widehat{f}_c(k,\tau) d\tau \right| (1+|k|^2)^{r/2}\right)\\
& \leq & \varepsilon^2 \sup_{t \in [0,T_0/\varepsilon^2]} \int_0^t \sup_{|k \pm k_c| \leq k_c/10} \left| e^{(\Lambda- \widetilde{\eta} \varepsilon^2) (t-\tau)} e^{-\widetilde{\eta} \varepsilon^2 \tau}\widehat{f}_c(k,\tau)(1+|k|^2)^{r/2} \right| d\tau
\\
& \leq & \varepsilon^2 \sup_{t \in [0,T_0/\varepsilon^2]} \int_0^t  C_{\Lambda} e^{(1-\widetilde{\eta} )\varepsilon^2 (t-\tau)} (1+|2 k_c|^2)^{r/2} \sup_{|k \pm k_c| \leq k_c/10} \left| e^{-\widetilde{\eta} \varepsilon^2 \tau} \widehat{f}_c(k,\tau) \right| d\tau
\\
& \leq &   C_2 \frac{1}{\widetilde{\eta}-1} \sup_{t \in [0,T_0/\varepsilon^2]}  \sup_{|k \pm k_c| \leq k_c/10}
\left|e^{-\widetilde{\eta} \varepsilon^2 t} \widehat{f}_c(k,t)\right| 
\\
& \leq &C_2 \frac{1}{\widetilde{\eta}-1} 
 \sup_{t \in [0,T_0/\varepsilon^2]} \|e^{-\widetilde{\eta} \varepsilon^2 t} \widehat{f}_c(\cdot,t)  \|_{L^{\infty}_{0} }
 \\
& = &C_2 \frac{1}{\widetilde{\eta}-1} \|  f_c \|_{
\mathcal{Y}^{0,\widetilde{\eta} \varepsilon^2}}
\end{eqnarray*}
for a constant $ C_2 > 0 $ independent of $ 0 < \varepsilon \ll 1 $ and $\widetilde{\eta} \geq 2$. 
\qed
\medskip

%
%

For $ \mathcal{K}_s $ we prove 
\begin{lemma} \label{lem73}
For all $ \widetilde{\eta} > 0 $, $ r \geq 2 $ there exist $ \varepsilon_0 > 0 $ and  $ C > 0 $ such that for all $ \varepsilon \in (0,\varepsilon_0) $ we have  
$$
\| \mathcal{K}_s f \|_{\mathcal{Y}^{r,\widetilde{\eta} \varepsilon^2} } \leq C  \|  f \|_{
\mathcal{Y}^{r-2,\widetilde{\eta} \varepsilon^2}}  .
$$
\end{lemma}
\noindent
{\bf Proof.} 
We set 
$$ 
W_s(\cdot,t) = e^{-\widetilde{\eta} \varepsilon^2 t} R_s(\cdot,t) 
$$
and find 
$$ 
\partial_t W_s(\cdot,t)  =  (\Lambda- \widetilde{\eta} \varepsilon^2) 
W_s(\cdot,t) 
+  e^{-\widetilde{\eta} \varepsilon^2 t}  (E_sf) (\cdot,t). 
$$ 
Applying the variation of constant formula yields
\begin{eqnarray*}
W_s(t) & = & \int_0^t e^{ (\Lambda- \widetilde{\eta} \varepsilon^2)  (t-\tau)} e^{-\widetilde{\eta} \varepsilon^2 \tau} (E_s f)(\tau) d\tau  .
\end{eqnarray*}
For the semigroup $ (e^{t\Lambda} )_{t \geq 0} $ applied 
to the $ E_s $-part we find 
$$ 
\| e^{(\Lambda- \widetilde{\eta} \varepsilon^2)t} W_{s,0} \|_{\mathcal{X}^r} \leq C_{\Lambda} e^{-\sigma t} \| W_{s,0}  \|_{\mathcal{X}^r}
$$ 
%
with a constant $ C_{\Lambda} = \mathcal{O}(1) $ and 
$ \sigma = \mathcal{O}(1) > 0 $.
After the diagonalization, 
the semigroup can be written as 
$$
e^{(\widehat{\Lambda}(k)- \widetilde{\eta} \varepsilon^2)t} = S^{-1}(k) \textrm{diag}(e^{\lambda_1(k) t},e^{\lambda_2(k) t}) S(k)
$$
where 
\begin{equation} \label{Sesti}
\sup_{|k \pm k_c| \geq k_c/10} \| S^{-1}(k) \|_{\C^2 \to \C^2} + \sup_{|k \pm k_c| \geq k_c/10} \| S(k) \|_{\C^2 \to \C^2} 
\leq C_S < \infty
\end{equation}
due to the asymptotics of $ \widehat{\Lambda}(k) $ for $ |k| \to  \infty $.
Then for $ \widehat{U}(k,t)  = S^{-1}(k)  \widehat{W}(k,t)   $ and 
$ \widehat{g}(k,t) = S^{-1}(k) \widehat{E}_s(k) \widehat{f}(k,t) $, 
we estimate for the components for each $ r \geq 2 $ that
\begin{eqnarray*}
\lefteqn{\sup_{t \in [0,T_0/\varepsilon^2]}\|(1-\chi_{|k| \geq k_0}(\cdot))\widehat{U}_j(\cdot,t)  \|_{L^{\infty}_r} }\\ & = & 
\sup_{t \in [0,T_0/\varepsilon^2]}\sup_{|k| \geq k_0}\left( \left| \int_0^t e^{\lambda_j(k)  (t-\tau)}e^{-\widetilde{\eta} \varepsilon^2 \tau}\widehat{g}_j(k,\tau) d\tau \right| (1+|k|^2)^{r/2}\right)\\
& = & \sup_{|k| \geq k_0} \sup_{t \in [0,T_0/\varepsilon^2]}\left( \left|\int_0^t e^{\lambda_j(k)  (t-\tau)} e^{-\widetilde{\eta} \varepsilon^2 \tau}\widehat{g}_j(k,\tau) d\tau \right| (1+|k|^2)^{r/2}\right)
\\
& \leq & \sup_{|k| \geq k_0} \sup_{t \in [0,T_0/\varepsilon^2]}\left( \int_0^t \left|e^{\lambda_j(k)  (t-\tau)} (1+|k|^2)\right| d\tau \right. \\&& 
\left. \qquad \qquad \times \sup_{\tau \in [0,T_0/\varepsilon^2]} \left(e^{-\widetilde{\eta} \varepsilon^2 \tau}|\widehat{g}_j(k,\tau) | (1+|k|^2)^{(r-2)/2}\right)\right)
\\
& \leq & \sup_{|k| \geq k_0} \sup_{t \in [0,T_0/\varepsilon^2]}\left( \frac{1+|k|^2}{\lambda_j(k)} \left.e^{\lambda(k)  (t-\tau)} \right|_{\tau=0}^t
\sup_{\tau \in [0,T_0/\varepsilon^2]} \left(e^{-\widetilde{\eta} \varepsilon^2 \tau}|\widehat{g}_j(k,\tau)  | (1+|k|^2)^{(r-2)/2}\right)\right)
\\
& \leq &C\sup_{|k| \geq k_0} \sup_{t \in [0,T_0/\varepsilon^2]}\left( 
 e^{-\widetilde{\eta} \varepsilon^2 t}|\widehat{g}_j(k,t) | (1+|k|^2)^{(r-2)/2}\right)
 \\
& =  & C\sup_{t \in [0,T_0/\varepsilon^2]} \sup_{k \in \R}\left( 
 e^{-\widetilde{\eta} \varepsilon^2 t}|\widehat{g}_j(k,t)  | (1+|k|^2)^{(r-2)/2}\right) \\
 & =  &C\sup_{t \in [0,T_0/\varepsilon^2]} \|e^{-\widetilde{\eta} \varepsilon^2 t}\widehat{g}_j(\cdot,t)  \|_{L^{\infty}_{r-2} },
\end{eqnarray*}
where we used $ \lambda_j(k) \leq \min(-\sigma,C-\alpha k^2) $
due to \eqref{eq:asymptotics_lambda}
for some $\alpha > 0$ depending only on the coefficients $a, b$ and $d$.
For $ |k| \leq k_0 $ the estimates are straightforward due to the compact support in Fourier space and the exponential decay of the semigroup in the $ E_s $-part.
\qed
\medskip

%
%

The error estimates will follow with a fixed point argument. From 
\eqref{erreq1}-\eqref{erreq2} we obtain
\begin{eqnarray} \label{erreq1inv}
R_c & = & F_c( R_c,R_s),
\\ \label{erreq2inv}
R_s & = &  F_s( R_c,R_s),
\end{eqnarray}
where
\begin{eqnarray*} 
 F_c( R_c,R_s) & = & \mathcal{K}_c( L_{2,c}(\varepsilon \psi,\varepsilon^2 R)
+ N_{2,c}(\varepsilon \psi,\varepsilon^2 R)  
 \\ && \qquad + \varepsilon^{-2} E_c N_3(\varepsilon \psi,\varepsilon^2 R) + \varepsilon^{-2} E_c \textrm{Res}(\varepsilon \psi)), 
 \nonumber \\
 F_s( R_c,R_s) & = & \mathcal{K}_s ( L_{2,s}(\varepsilon \psi,\varepsilon^2 R)
+ N_{2,s}(\varepsilon \psi,\varepsilon^2 R)   \\ && \qquad  + \varepsilon^{-3} E_s N_3(\varepsilon \psi,\varepsilon^2 R)
+ \varepsilon^{-3} E_s \textrm{Res}(\varepsilon \psi)).
\end{eqnarray*}
We prove that the mapping $ (R_c,R_s) \mapsto (F_c( R_c,R_s),F_s( R_c,R_s)) $
is a contraction in a ball in the space 
$$ 
\mathcal{Y}^{r,\widetilde{\eta} \varepsilon^2} \times \mathcal{Y}^{r,\widetilde{\eta} \varepsilon^2}
$$ 
equipped with norm
$$ 
(R_c,R_s) \mapsto   \| R_c \|_{\mathcal{Y}^{r,\widetilde{\eta} \varepsilon^2} } + 
\alpha \| R_s \|_{\mathcal{Y}^{r,\widetilde{\eta} \varepsilon^2} } 
$$ 
where $ \alpha > 0 $ is suitably chosen below to get rid of the Jordan block structure of $ 
(F_c,F_s) $.
With $  Z =  \| R_c \|_{\mathcal{Y}^{r,\widetilde{\eta} \varepsilon^2} } + 
\alpha \| R_s \|_{\mathcal{Y}^{r,\widetilde{\eta} \varepsilon^2} }  $ we estimate, using Lemma \ref{lem72} and \ref{lem73},
\begin{eqnarray*}  
 \| F_c( R_c,R_s)  \|_{\mathcal{Y}^{r,\widetilde{\eta} \varepsilon^2} } 
 & \leq & C (\widetilde{\eta}-1)^{-1} (
 \| R_c \|_{\mathcal{Y}^{r,\widetilde{\eta} \varepsilon^2} } + \| R_s \|_{\mathcal{Y}^{r,\widetilde{\eta} \varepsilon^2} }) \\ && + C (\widetilde{\eta}-1)^{-1} e^{2 \widetilde{\eta}} \varepsilon (\| R_c \|_{\mathcal{Y}^{r,\widetilde{\eta} \varepsilon^2} } + \| R_s \|_{\mathcal{Y}^{r,\widetilde{\eta} \varepsilon^2} })^2 
 \\ &&+ C (\widetilde{\eta}-1)^{-1}  C_{res}
 \\ & \leq & C (\widetilde{\eta}-1)^{-1} \alpha^{-1} Z + C e^{2 \widetilde{\eta}} \varepsilon ( \alpha^{-1} Z)^2 + C C_{res}
\end{eqnarray*}
and 
\begin{eqnarray*}  
\alpha \| F_s( R_c,R_s)  \|_{\mathcal{Y}^{r,\widetilde{\eta} \varepsilon^2} } 
 & \leq & C \alpha (
 \| R_c \|_{\mathcal{Y}^{r,\widetilde{\eta} \varepsilon^2} } + \varepsilon \| R_s \|_{\mathcal{Y}^{r,\widetilde{\eta} \varepsilon^2} }) \\ && \qquad + C e^{2 \widetilde{\eta}}  \varepsilon (\| R_c \|_{\mathcal{Y}^{r,\widetilde{\eta} \varepsilon^2} } + \| R_s \|_{\mathcal{Y}^{r,\widetilde{\eta} \varepsilon^2} })^2 + C C_{res}
\\ &  \leq & 
 C \alpha Z + C \varepsilon Z + C e^{2 \widetilde{\eta}}  \varepsilon Z^2 + C \alpha C_{res}. 
\end{eqnarray*}
In these estimates we used the fact that the approximations $ \psi_c $ and $ \psi_s $ 
are uniformly bounded on the time interval $ [0,T_0/\varepsilon^2] $ such 
that no factor $ e^{\widetilde{\eta}} $ occurs for the terms where $ R $ appears linearly.

Therefore, the right hand side maps a ball of radius $ r_0 $ in itself for 
$ \alpha> 0 $ chosen sufficiently 
small, then $   \widetilde{\eta}  $ chosen sufficiently large,  and finally $ \varepsilon > 0 $ chosen sufficiently 
small.
With the same argument  the right hand side is a contraction in this  ball of radius $ r_0 $.
Hence, the exists a unique fixed point in this ball for this mapping. 
\medskip

Therefore, we have established the following approximation result.
\begin{theorem}
For all $ C_{GL} > 0 $, $ T_0 > 0 $, $ r \geq 2 $ there exist $ \varepsilon_0 > 0 $ 
and $ C_2 > 0 $ such that the following holds.
Let $ A \in C([0,T_0],\mathcal{X}^{r+6}) $ be a solution of the 
Ginzburg-Landau equation \eqref{GL} satisfying 
\eqref{eq19}.
Then for all
$ \varepsilon \in (0,\varepsilon_0) $ there are solutions 
$ V $ of the Gray-Scott-Klausmeier system \eqref{OS}
such that 
$$ 
\sup_{t \in [0,T_0/\varepsilon^2]} \| V(\cdot,t) - \varepsilon \psi(\cdot,t) \|_{\mathcal{X}^{r}} \leq C_2 \varepsilon^2.
$$
\end{theorem}
Since we can estimate 
$$ 
\| V(\cdot,t) - \varepsilon \psi_{GL}(\cdot,t)  \|_{\mathcal{X}^{r}} \leq \| V(\cdot,t) - \varepsilon \psi(\cdot,t)  \|_{\mathcal{X}^{r}}
+ \|  \varepsilon \psi_{GL}(\cdot,t) - \varepsilon \psi(\cdot,t) \|_{\mathcal{X}^{r}}
$$ 
and since 
$$
\| R \|_{C_b^0} \leq \| \widehat{R} \|_{L^1} \leq 
C \| \widehat{R} \|_{L^{\infty}_r} = C  \| R \|_{\mathcal{X}^{r}}
$$ 
for $ r > 1 $ we also have 
\begin{theorem}
For all $ C_{GL} > 0 $, $ T_0 > 0 $, $ r > 1 $ there exist $ \varepsilon_0 > 0 $ 
and $ C_2 > 0 $ such that the following holds.
Let $ A \in C([0,T_0],\mathcal{X}^{r+6}) $ be a solution of the 
Ginzburg-Landau equation \eqref{GL} satisfying \eqref{eq19}. Then for all
$ \varepsilon \in (0,\varepsilon_0) $ there are solutions 
$ (v,w) $ of the Gray-Scott-Klausmeier system \eqref{gsk1}-\eqref{gsk2}
such that 
$$ 
\sup_{t \in [0,T_0/\varepsilon^2]} 
\sup_{x \in \R}
\| (v,w)(x,t) - ((v^\star,w^\star) + \varepsilon \psi_{GL}(x,t)) \|_{\R^2} \leq C_2 \varepsilon^2.
$$
\end{theorem}

\section{The general situation}
\label{secgeneral}

The previous approach works for general quasilinear pattern forming
reaction-diffusion-advection systems, too.
We consider 
\begin{equation}  \label{OS1} 
\partial_t U = D(U) \Delta U + f(U,\nabla U)
\end{equation}
with $ t \geq 0 $, $ x \in \R $, $ U(x,t) \in \R^m $, 
$ D(U) = \textrm{diag}(d_1(U),\ldots,d_m(U)) $ with $ d_j:\R^m \to \R $ 
smooth and $ d_j(U) > 0 $ for all $ j = 1,\ldots,m $, and finally 
smooth $ f: \R^m \times \R^m \to \R^m $.

Similar to Section \ref{sec2} we assume 
\medskip

{\bf (ASS2)} There exists a $ t $- and $ x $-independent fixed point $ U^* $. The linearization 
$$
\partial_t V = D(U^*) \Delta  V + \partial_1 f(U^*,0) V + 
\partial_2 f(U^*,0) \nabla V
$$ 
at $ U^* $ shows a Turing instability at a wave number $ k_c > 0 $, cf. Figure \ref{fig1}.  

Setting $ U = U^* + V $ gives  a system \eqref{OS} but now 
with additional terms $ \mathcal{O}(\| V \|^4) $.
Since terms of order $ \mathcal{O}(\| V \|^4) $ do not play a role in the derivation of the Ginzburg-Landau equation, Section \ref{sec3} 
transfer line for line from Gray-Scott-Klausmeier system 
\eqref{OS} to general 
reaction-diffusion-advection systems \eqref{OS1}.
There are no changes in Section \ref{sec4}. 
There are minor changes in Section \ref{sec5} due to the  
$ \mathcal{O}(\| V \|^4) $ terms. The improved ansatz is the same 
but the equations for $ A_{\pm 1,1} $ change slightly.
There are no changes in Section \ref{sec6} and also in Section \ref{sec7} if only polynomial $ f = f(U) $ is considered.
\begin{remark}{\rm
In case of polynomial $ f = f(U) $ the factor $ e^{\widetilde{\eta}} $ 
coming from Corollary \ref{collconvy} plays no role as we have seen.
However, for terms of order $ \mathcal{O}((\| R_c\|+\|R_s\|)^n) $ we would have a pre-factor 
$ e^{(n-1)\widetilde{\eta}} $. But this is not a problem since we have at least a pre-factor 
$ \varepsilon^{n-1} $, too, which still allows us to make general nonlinearities small. 
}\end{remark}
The estimate \eqref{Sesti} holds due to the fact that we assumed
that $ D(u) $ is a diagonal matrix.
Therefore, we have  
\begin{theorem}
Consider \eqref{OS1} and assume that the assumption 
{\bf (ASS2)} is fulfilled.
Then for all $ C_{GL} > 0 $, $ T_0 > 0 $, $ r > 2 $ there exist $ \varepsilon_0 > 0 $ 
and $ C_2 > 0 $ such that the following holds.
Let $ A \in C([0,T_0],\mathcal{X}^{r+6}) $ be a solution of the 
Ginzburg-Landau equation \eqref{GL} satisfying \eqref{eq19}. Then for all
$ \varepsilon \in (0,\varepsilon_0) $ there are solutions 
$ U $ of the 
general quasilinear pattern forming
reaction-diffusion-advection system
\eqref{OS1}
such that 
$$ 
\sup_{t \in [0,T_0/\varepsilon^2]} 
\sup_{x \in \R}
\| U(x,t) - (U^*+ \varepsilon \psi_{GL}(x,t)) \|_{\R^m} \leq C_2 \varepsilon^2.
$$
\end{theorem}
\begin{remark}{\rm
We assumed that the diffusion 
matrix $ D(U) $ is a diagonal matrix which would exclude 
cross-diffusion. In fact this can be weakened
by assuming that all eigenvalues of $ D(U^*) $ 
have negative real part. If all eigenvalues are different 
the estimate \eqref{Sesti} follows.}
\end{remark}
\begin{remark}{\rm
Following \cite{SU17book}
it is obvious that the present proof
can be transferred almost one-to-one
to the situation of a Turing-Hopf instability.
}
\end{remark}

\section{Discussion and Outlook}

\label{sec8}

The above proof that the Ginzburg-Landau equation approximates pattern-forming systems near the first instability differs only slightly in structure from the semilinear case.
Essentially, there are only three  but crucial points. 

The first point is the choice of the function space. In the semilinear case, for example,
Sobolev spaces
$ H^r $, spaces of continuously differentiable functions $ C^r_b $, or
uniformly locally integrable Sobolev spaces
$ H^r_{l,u} $ can be chosen. The function space we have chosen,
$ \mathcal{X}^r $, has the same disadvantage as Sobolev spaces, 
namely that the functions which it contains vanish at infinity.
This means that these spaces do not contain any spatially periodic functions or fronts.

The second point justifies the choice of our space.
In this space, an optimal regularity result applies, namely Lemma
\ref{lem73}.
In the semilinear case, an estimate such as
$ \| e^{\Lambda t} R_s \|_{H^{r+\theta}} \leq C \max(1,t^{-\theta/2}) 
\|  R_s \|_{H^{r}} $ would be used. If the nonlinearity loses two derivatives, then $ \theta = 1 $ is required, which leads to a non-integrable singularity in the variation of the constant formula. This estimate can be replaced in the spaces $ \mathcal{X}^r $ by the optimal regularity result Lemma
\ref{lem73}.

The third point concerns the final estimates following the proof of Lemma
\ref{lem73}.
In the semilinear case and also in \cite{Zi14PhD}, estimates are found here using Gronwall's inequality.
These are also used as a priori estimates to conclude the long-term existence of the solutions in connection with a local existence and uniqueness theorem. This would also be possible for the spaces we have chosen. However, we have decided to present an approach that can also be applied to more complicated situations,
and in particular can also be applied when continuity in time is replaced by Hölder continuity in time, cf. \cite{Lunardi}.
The fixed point argument presented works whenever  an
optimal regularity result in the sense of Lemma
\ref{lem73} is available.

The goal of future research will be to apply the present scheme to more complicated, pattern-forming, quasi-linear
systems, such as the B\'enard problem with velocity- and temperature-dependent viscosity. To work in function spaces whose elements do not vanish at infinity, Lemma \ref{lem73} will be replaced by 
more advanced optimal regularity results which we have been avoided 
here in order to have an easy to use approximation result avoiding
all these functional analytic tools.

\bibliographystyle{alpha}
\bibliography{GLbib}

\end{document}